\documentclass[preprint,11pt]{elsarticle}

\usepackage{amsmath}
\usepackage{algorithm}
\usepackage{amsmath}
\usepackage{amssymb}
\usepackage{graphics}
\usepackage{graphicx}
\usepackage[position = top]{subfig}
\usepackage{hyperref}
\usepackage{booktabs}

\usepackage{caption}
\newcommand{\bb}{\mathbf b}
\newcommand{\be}{\mathbf e}
 
\newcommand{\bg}{\mathbf g}
\newcommand{\bs}{\mathbf s}
\newcommand{\bu}{\mathbf u}

\newcommand{\bx}{\mathbf x}
\newcommand{\by}{\mathbf y}

\newcommand{\zero}{\mathbf 0}
\newcommand{\calo}{\mathcal O}
\newcommand{\eps}{\varepsilon}
\newcommand{\R}{\mathbb R}

\newcommand{\PP}{\mathbb P}

\newcommand{\wt}{\widetilde}
\newcommand{\ph}{\phantom}
\newcommand{\aux}{\omega} 

\newcommand{\mtxa}[2]{
\left[ \begin{array}{#1} #2 \end{array} \right]} 
\newcommand{\smtxa}[2]{
{\mbox{\scriptsize $\left[\!\! \begin{array}{#1} #2 \end{array} \!\! \right]$}}}
\DeclareMathOperator*{\argmin}{argmin}

\usepackage{amsthm}
\newtheorem{theorem}{Theorem}
\newtheorem*{acknowledgement}{Acknowledgement}

\newtheorem{example}[theorem]{Example}
\newtheorem{lemma}[theorem]{Lemma}
\newtheorem{proposition}{Proposition}

\usepackage{xcolor}
\newcommand{\gf}[1]{{\color{black}#1}}

\journal{-}

\begin{document}

\begin{frontmatter}

\title{A homogeneous Rayleigh quotient with applications in gradient methods}

\author[inst1]{Giulia Ferrandi}

\affiliation[inst1]{organization={Department of Mathematics and Computer Science, TU Eindhoven},
addressline={\\ PO Box 513}, 
city={Eindhoven},
postcode={5600MB}, 
country={The Netherlands}}

\author[inst1]{Michiel E. Hochstenbach}

\begin{abstract}
Given an approximate eigenvector, its (standard) Rayleigh quotient and harmonic Rayleigh quotient are two well-known approximations of the corresponding eigenvalue.
We propose a new type of Rayleigh quotient, the homogeneous Rayleigh quotient, and analyze its sensitivity with respect to perturbations in the eigenvector. Furthermore, we study the inverse of this homogeneous Rayleigh quotient as stepsize for the gradient method for unconstrained optimization.
The notion and basic properties are also extended to the generalized eigenvalue problem.
\end{abstract}

\begin{highlights}
\item We define a new homogeneous Rayleigh quotient as the minimizer of a residual for eigenvalue problems.
\item A nonlinear Galerkin condition for this homogeneous Rayleigh quotient is derived.
\item Asymptotic bounds on the relative error to an eigenvalue are obtained.
\item The quotient is compared with standard and harmonic Rayleigh quotients.
\item We study the inverse of this quotient as a stepsize for gradient methods.
\end{highlights}

\begin{keyword}
Homogeneous Rayleigh quotient \sep secant condition \sep eigenvalue problem \sep projective coordinates \sep unconstrained optimization \sep generalized eigenvalue problem
\MSC 65F15 \sep 65F10 \sep 65K05 \sep 90C20 \sep 90C30 \sep 65F50
\end{keyword}

\end{frontmatter}

\section{Introduction} \label{sec:intro}
Let $A$ be an $n \times n$ symmetric matrix with eigenvalues $\lambda_1 \le \cdots \le \lambda_n$. We are interested in the eigenproblem $A\bx = \lambda \bx$.
Let $\bu \in \R^n$ be an approximate eigenvector to $\bx$ with unit 2-norm.
Traditionally, the Rayleigh quotient
\begin{equation} \label{rq}
\theta = \frac{\bu^T\!A\bu}{\bu^T\bu} \qquad (= \bu^T\!A\bu)
\end{equation}
is the standard approach to determine an eigenvalue approximation corresponding to $\bu$, which yields the corresponding eigenvalue when $\bu$ is an exact eigenvector.
This $\theta$ satisfies two (related) optimality properties, a Galerkin condition (cf., e.g., \cite[pp.~13--14]{Par98})
\begin{equation} \label{opt1a}
A\bu - \theta \bu \perp \bu,
\end{equation}
and a minimum residual condition
\begin{equation} \label{opt1b}
\theta = \argmin_{\gamma} \, \|A\bu-\gamma \, \bu\|.
\end{equation}
A second well-known quantity is the harmonic Rayleigh quotient, which is, for instance, sometimes considered when one is interested in an eigenvalue near a given target $\tau \in \R$. This is defined as (cf., e.g., \cite[p.~294]{Ste01})
\begin{equation} \label{hrq-tau}
\wt \theta_\tau = \frac{\bu^T(A-\tau I)\,A\bu}{\bu^T(A-\tau I)\,\bu},
\end{equation}
provided that the denominator is nonzero.
This quotient is also equal to the eigenvalue when $\bu$ is an eigenvector.
It satisfies two optimality conditions: a Galerkin condition $(A-\wt \theta_{\tau} I)\, \bu \perp (A-\tau I)\,\bu$, and a minimum residual condition (or secant condition) $\wt \theta_{\tau} = \argmin_{\gamma} \|(\gamma-\tau)^{-1} (A-\tau I)\,\bu-\bu\|$ (cf.~\cite[Prop.~5]{FHK23}).

In this paper, we mainly focus on the case $\tau = 0$; for nonzero $\bu^T\!A\bu$, the harmonic Rayleigh quotient is given by (cf.~\cite{morgan1991computing})
\begin{equation} \label{hrq}
\wt \theta = \frac{\bu^T\!A^2\,\bu}{\bu^T\!A\bu}.
\end{equation}
For $\tau = 0$, the Galerkin condition becomes 
\begin{equation} \label{opt2a}
A\bu - \wt \theta \bu \perp A\bu,
\end{equation}
and the minimum residual condition is
\begin{equation} \label{opt2b}
\wt \theta = \argmin_{\gamma} \, \|\gamma^{-1} A\bu-\bu\|.
\end{equation}
Condition \eqref{opt2b} is easily verified by setting the derivative of $\|\gamma^{-1} A\bu-\bu\|$ with respect to $\gamma$ to zero.
Property \eqref{opt2b} is perhaps less widely known.
Both Properties~\eqref{opt1b} and \eqref{opt2b} have been exploited as secant conditions to determine a stepsize in gradient-type optimization methods; for instance for the well-known Barzilai and Borwein steplengths \cite{bb1988}. We return to this topic in Section~\ref{sec:step}; see also V\"omel \cite{vomel2010note}.

As an alternative to \eqref{opt1b} and \eqref{opt2b}, we propose and study a new {\em homogeneous} type of Rayleigh quotient in this paper.
\gf{We work in the projective space $\PP = \PP_1(\R)$, also called the \emph{projective line}. This is a quotient space of $\mathbb R^2\setminus\{\bf 0\}$, with the following equivalence relation: two points $[\alpha_1,\alpha_2]^T$, $[\wt\alpha_1,\wt\alpha_2]^T\in \mathbb R^2\setminus\{\bf 0\}$ are equivalent if there exists a nonzero $\zeta$ such that $[\alpha_1,\alpha_2]^T = \zeta \, [\wt\alpha_1,\wt\alpha_2]^T$. The coordinates of the representative of a class $(\alpha_1,\alpha_2)$ are called \emph{homogeneous coordinates}.
The pair $(1,0)$ corresponds to the point at infinity. All other equivalence classes $(\alpha_2\ne 0)$ can be identified with a real number $\alpha = \alpha_1/\alpha_2$.}

Homogeneous techniques for eigenvalue problems have already been developed by various authors.
Stewart and Sun \cite[p.~283]{SSu90} exploit homogeneous coordinates for the generalized eigenvalue problem $A\bx = \lambda B\,\bx$, for $B \in \R^{n \times n}$, since this is ``especially convenient for treating infinite eigenvalues'' \cite{Ste75}.
While the standard eigenvalue problem $A\bx = \lambda \bx$ does not involve infinite eigenvalues, the harmonic Rayleigh quotient \eqref{hrq} may be infinite.
Dedieu and Tisseur \cite{Ded97,DTi03} exploit homogeneous techniques to study perturbation theory for generalized and polynomial eigenproblems.
Inspired by this, a homogeneous Jacobi--Davidson method for subspace expansion has been proposed in \cite{HNo08}; this may be viewed as an inexact Newton method.

The contributions and outline of this paper are the following. We define a new \emph{homogeneous Rayleigh quotient} as the minimizer of the residual quantity
\[\min_{(\alpha_1,\,\alpha_2) \, \in \, \mathbb P} \ \frac{\|\alpha_1 \, \bu - \alpha_2 \, A\bu\|}{\gf{\sqrt{\alpha_1^2+\alpha_2^2}}}.\]
In Section~\ref{sec:homo}, we show that there is a closed-form solution to express this homogeneous Rayleigh quotient, and that it can also be obtained from a certain {\em nonlinear} Galerkin condition.
To quantify the quality of the homogeneous Rayleigh quotient as an approximate eigenvalue, a bound for the (chordal) distance between the homogeneous Rayleigh quotient and the spectrum of $A$ is derived. 

We compare the homogeneous Rayleigh quotient with the standard and harmonic Rayleigh quotients. Given the relation between the homogeneous residual quantity and the two residual quantities \eqref{opt1b} and \eqref{opt2b} (see also Section~\ref{sec:relation}), the homogeneous Rayleigh quotient may be seen as the mediator between the two other Rayleigh quotients. In Section~\ref{sec:sens}, we highlight this fact in the theoretical and experimental study on the relative errors of the three quotients. Asymptotic bounds on the accuracy of the three quotients are also derived and compared. Interestingly, the homogeneous Rayleigh quotient may be more accurate than its two competitors in some situations. 

In Section~\ref{sec:gep}, we extend the definition of the homogeneous Rayleigh quotient to the generalized eigenvalue problem $(A, B)$, with $A$ symmetric and $B$ symmetric positive definite (SPD).
Finally, in Section~\ref{sec:step}, we propose a new stepsize for gradient methods based on the homogeneous Rayleigh quotient. Conclusions are drawn in Section~\ref{sec:con}.

\gf{We would like to stress that the stepsize \eqref{eq:homobb} presented in Section~\ref{sec:step} has been recently proposed first by Li, Zhang, and Xia \cite{LZX22} (with preliminary version \cite{LX21} in 2021).
The results in our Section~\ref{sec:step} have been developed independently and the  derivation is different: whereas the approach by Li, Zhang, and Xia is via a total least squares technique, ours is based on exploiting homogeneous coordinates.
We will discuss some more details in Section~\ref{sec:step}.}

\section{A homogeneous Rayleigh quotient} \label{sec:homo}
For a given approximate eigenvector $\bu$, we consider the minimization of the \emph{homogeneous} residual (as mediator of the two residual quantities in \eqref{opt1b} and \eqref{opt2b}):
\begin{equation} \label{homo-sec}
\min_{(\alpha_1,\,\alpha_2) \, \in \, \mathbb P} \ \frac{\|\alpha_1 \, \bu - \alpha_2 \, A\bu\|}{\gf{\sqrt{\alpha_1^2+\alpha_2^2}}} = \min_{\alpha_1^2+\alpha_2^2=1} \ \|[\bu \ \ - \!\!A\bu] \smtxa{c}{\alpha_1 \\ \alpha_2}\|.
\end{equation}
\gf{Note that the expression on the left-hand side is well defined on $\mathbb P$. The equality between the two optimization problems can be seen by taking $(\alpha_1,\alpha_2) \leftarrow (\alpha_1^2 + \alpha_2^2)^{-1/2} \ (\alpha_1, \alpha_2)$ as representative of the equivalence class. Then we minimize the objective function on $\{[\alpha_1,\alpha_2]^T\in\mathbb R^2\mid \alpha_1^2 + \alpha_2^2 = 1\}$.}
For nonzero $\alpha_2$ we call the ratio between the coordinates of the solution $\alpha = \alpha_1/\alpha_2$ the \emph{homogeneous Rayleigh quotient}.

\subsection{Key properties of the homogeneous Rayleigh quotient} We discuss some properties of the solution to \eqref{homo-sec}. First, let us introduce the $n \times 2$ matrix $C = [\bu \ \ \, -\!\!A\bu]$ and the associated $2 \times 2$ matrix
\begin{equation} \label{m22}
C^T C = \mtxa{cc}{\ph-\bu^T\bu & -\bu^T\!A\bu \\ -\bu^T\!A\bu & \ph-\bu^T\!A^2\bu} = \bu^T\!A\bu \cdot \mtxa{cc}{\ph-\theta^{-1} & -1 \\ -1 & \ph-\wt \theta},
\end{equation}
where the second equality holds when $\bu^T\!A\bu$ is nonzero.
Although we usually impose $\|\bu\| = 1$, we will not exploit this simplification in view of more general use in Section~\ref{sec:step}. 

In the next proposition, assuming that $\bu^T\!A\bu \ne 0$, we can provide an explicit formula for the homogeneous Rayleigh quotient $\alpha$. In addition, we show that $\alpha$ is an eigenvalue of $A$ when $\bu$ is an eigenvector, and that it is located between the Rayleigh quotient and the harmonic Rayleigh quotient, depending on the sign of $\bu^T\!A\bu$.

\begin{proposition} \label{prop:prop}
Suppose that $\bu^T\!A\bu \ne 0$ and denote the smallest eigenvalue of $C^TC$ by $\mu$. Then the following properties hold.
\begin{itemize}
\item[(i)] The value $\mu$ is a simple eigenvalue of $C^TC$, or, equivalently, $\sqrt{\mu}$ is a simple singular value of $C$; its corresponding eigenvector $[\alpha_1, \ \alpha_2]^T$ is the unique minimizer of \eqref{homo-sec}. 
\item[(ii)] The value of $\mu$ is
\begin{align}
\mu = \tfrac12 \, \big[\bu^T\bu+\bu^T\!A^2\,\bu - \sqrt{\smash[b]{(\bu^T\bu-\bu^T\!A^2\,\bu)^2+4\,(\bu^T\!A\bu)^2}} \ \big]. \label{mu}
\end{align}
\item[(iii)] For the homogeneous Rayleigh quotient $\alpha$ we have
\begin{align} \label{eq:homorq}
\alpha &= \frac{\alpha_1}{\alpha_2} = \frac{\bu^T\!A\bu}{\bu^T\bu - \mu} = \frac{\bu^T\!A^2\,\bu - \mu}{\bu^T\!A\bu}\nonumber\\
&= \frac1{2\ \bu^T\!A\bu} \ \big[ \,\bu^T\!A^2\,\bu - \bu^T\bu + \sqrt{\smash[b]{(\bu^T\!A^2\,\bu - \bu^T\bu)^2 + 4\,(\bu^T\!A\bu)^2}} \,\big].
\end{align}
In terms of the standard and harmonic Rayleigh quotient,
\[
\alpha = \tfrac12 \, \big[ \, \wt \theta-\theta^{-1} + \sqrt{\smash[b]{(\wt \theta-\theta^{-1})^2 + 4}} \, \big]
= \big( \tfrac12 \, \big[ \, \theta^{-1}-\wt \theta + \sqrt{\smash[b]{(\theta^{-1}-\wt \theta)^2 + 4}} \, \big] \big)^{-1}.
\]
\item[(iv)] $\mu=0$ if and only if $\bu$ is an eigenvector of $A$. If $\bu$ is an eigenvector, then the corresponding homogeneous Rayleigh quotient $\alpha$ is the corresponding eigenvalue.
\item[(v)] $0 \le \mu < \min(\bu^T\bu, \ \bu^T\!A^2\bu)$.
\item[(vi)] We have the following bounds:
\begin{equation*}
\begin{cases}
\theta \le \alpha \le \wt \theta &\text{if}\quad\bu^T\!A\bu > 0, \\[0.5mm]
\wt \theta \le \alpha \le \theta &\text{if}\quad\bu^T\!A\bu < 0.
\end{cases}
\end{equation*}
\end{itemize}
\end{proposition}

\begin{proof}
It is well known that the minimizer of \eqref{homo-sec} is the eigenvector corresponding to the smallest eigenvalue of the matrix $C^T C$; or, equivalently, to the right singular vector of $C$ corresponding to the smallest singular value.
We start with the eigenvalue problem
\[C^TC\,\smtxa{c}{\alpha_1 \\ \alpha_2} = \mu \smtxa{c}{\alpha_1 \\ \alpha_2}.\]
The eigenvalue $\mu$ of \eqref{m22} is simple because $\bu^T\!A\bu$ is assumed nonzero, so that the discriminant of the quadratic characteristic polynomial is positive. In particular, this implies that the smallest eigenvector $[\alpha_1, \ \alpha_2]^T$ is well defined. The result in part (ii) is an explicit expression for the eigenvalue $\mu$. Its corresponding eigenvector $[\alpha_1, \ \alpha_2]^T$ is used in item (iii) to give an expression \gf{of}
the ratio $\alpha = \alpha_1/\alpha_2$, which gives the homogeneous Rayleigh quotient, that can also be written in terms of the standard and the harmonic Rayleigh quotients.

If $\mu = 0$, this means that $C$ has rank 1. This happens precisely when $\bu$ is an eigenvector. Moreover, if $\bu$ is an eigenvector corresponding to some eigenvalue $\lambda$, equations \eqref{mu} and \eqref{eq:homorq} give $\mu = 0$ and $\alpha = \lambda$ respectively. This concludes the proof of part (iv).

The nonnegativity of $\mu$ directly follows from the fact that $C^TC$ is positive semidefinite. Since $\bu^T\!A\bu \ne 0$ we have $\sqrt{\smash[b]{(\bu^T\bu-\bu^T\!A^2\,\bu)^2+4\,(\bu^T\!A\bu)^2}} > \vert \bu^T\bu-\bu^T\!A^2\,\bu \vert$, and therefore $\bu^T\bu -\mu,\ \bu^T\!A^2\bu -\mu > 0$.
It follows that $0\le \mu < \min(\bu^T\bu, \, \bu^T\!A^2\bu)$. In particular, this shows that the sign of $\alpha$ only depends on the quantity $\bu^T\!A\bu$. Given this fact, property (vi) is a direct consequence of part (v).
\end{proof}

In the exceptional case of $\bu^T\!A\bu = 0$, the matrix $C^TC$ is diagonal, with zero Rayleigh quotient and infinite harmonic Rayleigh quotient.
The homogeneous Rayleigh quotient $\alpha = (\alpha_1, \alpha_2)$ is either zero (that is, $(0,1)$), or the point at infinity $(1,0)$.
The situation $\bu^T\!A\bu = 0$ cannot happen if $A$ is definite, such as for quadratic optimization problems with an SPD Hessian; see also Section~\ref{sec:step}.
Moreover, $\bu^T\!A\bu = 0$ may not be a very common situation in the context of eigenproblems; this usually means that we are approximating a zero eigenvalue or that $\bu$ is a poor approximate eigenvector.

Rayleigh quotients have some natural invariance properties. It is easy to check that the standard Rayleigh quotient is invariant under nonzero scalings of $\bu$; it is invariant under nonzero scalings of $A$, and it shifts naturally with shifts of $A$.
The harmonic Rayleigh quotient (with target $\tau$) only satisfies the first property, but modified properties still hold, as follows.

\begin{proposition} \label{prop:hrq}
For any $\tau \in \mathbb{R}$, the following properties hold:
\begin{itemize}
\item[(i)] $\wt \theta_{\tau}$ is invariant under scaling of $\bu$: $\wt \theta_{\tau}(\zeta \, \bu) = \wt \theta_{\tau}(\bu)$ for $\zeta \ne 0$.
\item[(ii)] For $\wt \theta_{\tau}$ as function of scaling of $A$ we have
$\wt \theta_{\tau}(\zeta \, A) = \zeta \, \wt \theta_{\zeta^{-1} \, \tau}(A)$ for $\zeta \ne 0$.
\item[(iii)] For $\wt \theta_{\tau}$ as function of shifts of $A$ we have
$\wt \theta_{\tau}(A-\eta I) = \wt \theta_{\tau+\eta}(A) - \eta$ for $\zeta \in \R$.
\end{itemize}
\end{proposition}
\begin{proof}
This is easy to check from \eqref{hrq-tau}.
\end{proof}

As for the homogeneous Rayleigh quotient, it is easy to see (for instance, from \eqref{eq:homorq}) that the solution to the minimization problem \eqref{homo-sec} is invariant under nonzero scaling of $\bu$: as function of $\bu$, we have $\alpha(\zeta \, \bu) = \alpha(\bu)$ for $\zeta \ne 0$.
Second, while the standard and harmonic Rayleigh quotients are invariant under multiplications $A \to \zeta A$, there is no easy relation between $\alpha(\zeta A)$ and $\alpha(A)$. However, in view of Proposition~\ref{prop:prop}(vi), we establish these bounds:
\begin{equation*}
\begin{cases}
\zeta \, \theta(A) \le \alpha(\zeta A) \le \zeta \, \wt \theta(A) &\text{if}\quad\bu^T\!A\bu > 0, \\[0.5mm]
\zeta \, \wt \theta(A) \le \alpha(\zeta A) \le \zeta \, \theta(A) &\text{if}\quad\bu^T\!A\bu < 0.
\end{cases}
\end{equation*}
The same upper and lower bounds hold for $\zeta\alpha(A)$. Therefore, if the interval $[\zeta\theta(A),\, \zeta\wt \theta(A)]$ (or $[\wt\zeta\theta(A),\, \zeta\theta(A)]$) is small, we conclude that $\alpha(\zeta A) \approx \zeta \alpha(A)$. 

On the size of this interval, we note that when $\theta \le \wt \theta$ (which for instance holds in the SPD case), the ``relative size'' of the interval $[\theta,\, \wt \theta]$ can be expressed in terms of the angle between $\bu$ and $A\bu$, since 
$(\wt \theta-\theta) \, / \, \theta = \cos^{-2}(A\bu, \bu)-1 = \tan^2(A\bu, \bu)$.
This means that if $\angle(A\bu, \bu)$ is small, the interval is small, and all three Rayleigh quotients are close.
Still, Example~\ref{ex:smallex} will show that one quotient can still be considerably more accurate than the other two.

Moreover, we have
\[
\frac{\alpha - \theta}{\theta} = \tfrac12 \, \big[ \, 1 + \tan^2(A\bu, \bu)-\theta^{-2} + \sqrt{\smash[b]{(1 + \tan^2(A\bu, \bu)-\theta^{-2})^2 + 4\,\theta^{-2}}} \, - 2\big].
\]
For $\theta^{-1} \ll \tan(A\bu,\bu)$ (that is, for large eigenvalues), this expression is close to $\tan^2(A\bu, \bu)$, which means that $\alpha \approx \wt \theta$.
In contrast, if $\theta^{-1} \gg \max(\tan(A\bu,\bu), \, 1)$ (i.e., for small eigenvalues), one may check that $\frac{\alpha-\theta}{\theta} = \calo(\theta^2) \approx 0$, which implies that $\alpha \approx \theta$.
In conclusion, for a small $\angle(A\bu, \bu)$, the homogeneous Rayleigh quotient is close to the Rayleigh  quotient for small eigenvalues, and close to the harmonic Rayleigh quotient for large ones.
Similar remarks can be made when $\wt \theta < \theta$.

Simple manipulations of Proposition~\ref{prop:prop}(vi) also lead to some relations between the relative errors of the three Rayleigh quotients, when the estimated eigenvalues are either $\lambda_1$ or $\lambda_n$. If $A$ is SPD, then both the standard and harmonic Rayleigh quotients lie in the spectrum of $A$, i.e., $\theta,\,\wt\theta \in [\lambda_1,\lambda_n]$. Thus the following inequalities hold: 
\begin{equation}
\label{eq:sens-spd}
\theta - \lambda_1 \le \alpha - \lambda_1 \le \wt\theta - \lambda_1, \qquad \lambda_n -\wt\theta \le \lambda_n -\alpha \le \lambda_n -\theta.
\end{equation}
We conclude that, in the SPD case, the Rayleigh quotient is more accurate when we estimate the smallest eigenvalue of $A$, while the harmonic Rayleigh quotient is more accurate for the largest eigenvalue. The homogeneous Rayleigh quotient lies in between.
This may suggest the perhaps not so well-known result that to approximate large eigenvalues, the harmonic Rayleigh quotient may be preferable over the standard Rayleigh quotient.

For indefinite matrices, while $\theta\in [\lambda_1,\lambda_n]$, the harmonic and the homogeneous Rayleigh quotients can lie outside the spectrum. In particular, \eqref{eq:sens-spd} does not hold if $A$ is not SPD. Example~\ref{ex:smallex} shows that the harmonic and homogeneous Rayleigh quotients may overestimate the largest eigenvalues, and the sensitivity of the homogeneous Rayleigh quotient may be smaller than the other two. 

\begin{example}
\label{ex:smallex}
\rm Let $A = \text{diag}(-\frac23, \ \frac13, \ 2)$, and consider the eigenvector $\bx = [0, \, 0, \, 1]^T$ corresponding to $\lambda_3 = 2$. Given $\bu = [-0.02, \ 0.01, \ 1]^T$ as approximation to $\bx$, with $\angle(\bu, \bx) \approx 2.2 \cdot 10^{-2}$, we have that
\[
|\theta(\bu) - \lambda_3| \approx 1.2 \cdot 10^{-3}, \qquad
|\wt \theta(\bu) - \lambda_3| \approx 3.3 \cdot 10^{-4}, \qquad
|\alpha(\bu) - \lambda_3| \approx 1.6 \cdot 10^{-5}.
\]
The quotients are $\theta(\bu) \approx 1.9988 < \lambda_3$, $\wt \theta(\bu) \approx 2.0003 > \lambda_3$ and $\alpha(\bu) \approx 2.00002 > \lambda_3$, so the harmonic and the homogeneous Rayleigh quotients overestimate $\lambda_3$. Nevertheless, the homogeneous Rayleigh quotient is an order of magnitude more accurate than the harmonic Rayleigh quotient, and two orders more accurate than the standard Rayleigh quotient.
This example highlights that the homogeneous Rayleigh quotient might be more accurate, especially for exterior eigenvalues of indefinite matrices. We will show another example in Figure~\ref{fig:sigmas} in Section~\ref{sec:sens}.
\end{example}

When the wanted eigenvalue is in the interior of the spectrum, or if $A$ is indefinite, it is more difficult to derive general conclusions about the ordering of the relative errors. We will discuss this in more detail in Section~\ref{sec:sens}.

\subsection{A Galerkin condition for the homogeneous Rayleigh quotient}
In view of \eqref{opt1a} and \eqref{opt2a}, the question arises whether there exists a Galerkin (orthogonality) condition based on the span of $\bu$ and $A\bu$ for the homogeneous approach. Let us first express the homogeneous Rayleigh quotient as a solution to a quadratic equation. As a direct consequence, this property connects the homogeneous Rayleigh quotient and the harmonic Rayleigh quotient with target. 
\begin{proposition}\label{prop:quadeq}
Let $\bu^T\!A\bu \ne 0$. The homogeneous Rayleigh quotient \eqref{eq:homorq} solves
\begin{equation}
\label{eq:homopoly}
(\bu^T\!A\bu) \, \alpha^2 + (\bu^T\bu-\bu^T\!A^2\bu) \, \alpha - (\bu^T\!A\bu) = 0
\end{equation}
and satisfies $\alpha\ \bu^T\!A\bu > 0$. In addition, the homogeneous Rayleigh quotient is a harmonic Rayleigh quotient with target \eqref{hrq-tau}, where $\tau = -\alpha^{-1}$.
\end{proposition}
\begin{proof}
It is straightforward to check that \eqref{eq:homorq} is a solution to \eqref{eq:homopoly}. The solutions to \eqref{eq:homopoly} have opposite signs since their product is $-1$. Given that $0\le \mu < \min(\bu^T\bu, \, \bu^T\!A^2\bu)$ from Proposition~\ref{prop:prop}, the homogeneous Rayleigh quotient has the same sign as $\bu^T\!A\bu$, thus it corresponds to the solution to \eqref{eq:homopoly} for which $\alpha\ \bu^T\!A\bu > 0$.
The second part follows immediately from the fact that imposing $\wt\theta_{-\alpha^{-1}} = \alpha$ is equivalent to solving \eqref{eq:homopoly}.
\end{proof}

We remark that the second solution of \eqref{eq:homopoly} is related to the smallest eigenvector of $C^TC$, which maximizes (rather than minimizes) the homogeneous residual quantity \eqref{homo-sec}.

Another equivalent point of view on this proposition is the following.
Equation \eqref{eq:homopoly} can also be expressed as
\begin{equation}
\label{eq:homogal}
\bu^T(A + \alpha^{-1} I)\,(A-\alpha I)\,\bu = 0.
\end{equation}
Therefore, we have the {\em nonlinear Galerkin condition} $(A-\alpha I)\,\bu \perp (A+\alpha^{-1} I)\, \bu$, or, in homogeneous coordinates, 
\[
(\alpha_2 A-\alpha_1 I) \, \bu \perp (\alpha_1 A+ \alpha_2 I) \, \bu.
\]
This again highlights that $\alpha$ can be viewed as a harmonic Rayleigh quotient with $-\alpha^{-1}$ as target.

\subsection{Properties of the homogeneous residual quantity}
\label{sec:relation}
There is a connection between the two minimal residual conditions of the standard and harmonic Rayleigh quotients (\eqref{opt1b} and \eqref{opt2b}) and the quadratic equation \eqref{eq:homopoly}. 
The stationarity conditions for \eqref{opt1b} and \eqref{opt2b}, respectively, can be stated as
\[
\bu^T\bu - (\bu^T\!A\bu) \, \gamma^{-1} = 0, \qquad (\bu^T\!A\bu) \, \gamma - \bu^T\!A^2\bu = 0.
\]
This means that \eqref{eq:homogal} is a linear combination of the two stationarity conditions of the standard and the harmonic Rayleigh quotients. This fact, together with the bounds derived in Proposition~\ref{prop:prop}(vi), lets us interpret the homogeneous Rayleigh quotient as a ``mediator'' between the standard and the harmonic Rayleigh quotient.
Experiments in Section~\ref{sec:sensexp} will also highlight this mediating behavior.

It is an open question whether it is possible to relate the objective function \eqref{homo-sec} to a combination of the two residuals associated with the standard and the harmonic Rayleigh quotients. Nevertheless, it is possible to show that, when $\bu^T\!A\bu \ne 0$, the homogeneous residual quantity is smaller than the other two residual quantities (see \eqref{opt1b} and \eqref{opt2b}). \gf{As remarked after Proposition~\ref{prop:prop}, this hypothesis excludes the points $(0,1)$ and $(1,0)$ from the admissible minimizers. Therefore, \eqref{homo-sec} is restricted to $(\alpha_1,\alpha_2)\in\PP$, with $\alpha_1,\alpha_2\ne 0$. Then we use either $(\gamma, 1)$, or $(1, \gamma)$, for a nonzero $\gamma$, as class representatives. We derive the following minima, which are equivalent to \eqref{homo-sec} under $\bu^T\!A\bu \ne 0$:
\begin{equation}\label{eq:homores-alternative}
    \min_{\gamma\ne 0} \ \frac{1}{\sqrt{\gamma^2+1}} \, \|A\bu - \gamma \bu\| = \min_{\gamma\ne 0} \ \frac{\vert\gamma\vert}{\sqrt{\gamma^2+1}} \, \|\gamma^{-1}A\bu - \bu\|.
\end{equation}}
This means that the homogeneous residual quantity corresponds to the standard and harmonic residual quantities multiplied by a factor smaller than one. 

From \eqref{eq:homores-alternative}, we may also relate the homogeneous residual to a measure of distance between the homogeneous Rayleigh quotient and the spectrum of $A$.  
As stated in \cite[Thm.~4.5.1]{Par98} and \cite[Thm.~3.7]{vomel2010note}, for the standard Rayleigh quotient and the harmonic Rayleigh quotient there exists certain eigenvalues of $A$, say $\lambda$ and $\wt\lambda$, such that
\[
\vert\lambda - \theta\vert \le \frac{\|(A-\theta\, I)\,\bu\|}{\|\bu\|},\qquad \vert\wt\lambda^{-1} - \wt\theta^{-1}\vert \le \frac{\|(\wt\theta^{-1}I- A)\,\bu\|}{\|A\bu\|},
\]
provided that $A\bu \ne \zero$ in the second case.
These inequalities hold for any $\theta$ and $\wt\theta$, but the Rayleigh quotients $\theta$ and $\wt\theta$ have the advantage of minimizing the corresponding residual quantities.
A similar result can be stated for the homogeneous Rayleigh quotient in terms of the chordal metric (see, e.g., \cite[Ch.~6, Def.~1.20]{SSu90}). By applying \cite[Thm.~4.5.1]{Par98}, we get the following result:
for any $\alpha$, there exists an eigenvalue $\lambda$ of $A$ such that
\begin{equation} \label{eq:chord}
\frac{\vert\lambda - \alpha\vert}{\sqrt{1+\lambda^2} \ \sqrt{1+\alpha^2}} \le \frac{1}{\sqrt{1+\lambda^2} \ \|\bu\|}\cdot\frac{\|(A-\alpha\, I)\,\bu\|}{\sqrt{1+\alpha^2}},
\end{equation}
where the upper bound is minimized by the choice of the homogeneous Rayleigh quotient for $\alpha$.

Note that, as explained in \cite[Ch.~6]{SSu90}, the chordal metric behaves counterintuitively for large eigenvalues: while the chordal metric is small, the relative error can still be large.
For a fair comparison of the accuracy of the three Rayleigh quotients, we consider asymptotic bounds on their relative error when the vector $\bu$ approximates an eigenvector in the next section.

\section{Sensitivity analysis} \label{sec:sens}
We now study the sensitivity of the various Rayleigh quotients with respect to perturbations in the approximate eigenvector $\bu$.
Note that this sensitivity is different from that expressed by the condition number $\kappa(\lambda)$ of an eigenvalue, which is related to its perturbation as \gf{a} function of changes in the matrix $A$. We note that for a simple eigenvalue $\lambda$ of a symmetric $A$, it holds that $\kappa(\lambda) = 1$; this means that the eigenvalue is perfectly conditioned (see, e.g., \cite[p.~16]{Par98}).

Studying the sensitivity with respect to the approximate eigenvector for the standard and harmonic Rayleigh quotient is certainly not new (cf., e.g., \cite{sleijpen2003use,sleijpen2003optimal} and the references therein), but we derive a new result for the homogeneous Rayleigh quotient, and obtain expressions that allow an easy comparison of the three Rayleigh quotients. We will also comment on the differences between existing results and ours.

\subsection{Rayleigh quotient and harmonic Rayleigh quotient}
As ansatz, suppose $\bu=\bx+\be$ is an approximate eigenvector, where $\|\bx\| = 1$, $\be \perp \bx$, and $\eps := \|\be\|$ is small.
Then the perturbation of the standard Rayleigh quotient is
\begin{equation} \label{rq-sens}
\theta(\bu) = \frac{\bu^T\!A\bu}{\bu^T\bu} = \frac{\lambda + \be^T\!A\be}{1+\be^T\be}
= \lambda \, (1 + \be^T (\lambda^{-1} A-I)\,\be) + \calo(\eps^4),
\end{equation}
where we use the fact that $(1+t)^{-1} = 1 - t + \calo(t^2)$ for $t \to 0$.
For the harmonic Rayleigh quotient we have
\begin{align}
\wt \theta(\bu) = \frac{\bu^T\!A^2\,\bu}{\bu^T\!A\bu} & = \frac{\lambda^2+\be^T\!A^2\be}{\lambda+\be^T\!A\be}
= \lambda \, (1+\lambda^{-2} \, \be^T\!A^2\be) \, (1-\lambda^{-1} \, \be^T\!A\be) + \calo(\eps^4) \nonumber \\[1mm]
& = \lambda \, (1 + \be^T \lambda^{-1}A \, (\lambda^{-1} A-I)\,\be) + \calo(\eps^4). \label{hrq-sens}
\end{align}
For the harmonic Rayleigh quotient with target $\tau$ the expression becomes
\begin{equation} \label{hrq-tau-sens}
\wt \theta_\tau(\bu) = \lambda \, (1 + \be^T \, (\lambda-\tau)^{-1} (A-\tau I) \, (\lambda^{-1} A-I)\,\be) + \calo(\eps^4).
\end{equation}
In particular, we note that $\theta(\bx) = \wt \theta(\bx) = \wt \theta_\tau(\bx) = \lambda$. We also point out that we have studied similar expressions for {\em inverse} Rayleigh quotients (as stepsizes for gradient methods) in \cite{FHK23}.

To derive the lower and upper bounds on the sensitivity of the approximate eigenvalues, we make use of this standard result for symmetric operators.

\begin{lemma}
\label{lemma}
Suppose $\bu = \bx + \be$ is an approximate eigenvector corresponding to a simple eigenvalue $\lambda$, with $\|\bx\|=1$, $\be \perp \bx$, and $\eps = \|\be\|$, of a symmetric $A$. Let $p$ be a polynomial. Then
\[
\eps^2\,\min_{\lambda_i \ne \lambda} \ |p(\lambda_i)| \ \le \ |\be^T p(A) \, \be| \ \le \ \eps^2 \,\max_{\lambda_i \ne \lambda} \ |p(\lambda_i)|.
\]
\end{lemma}
\begin{proof}
Since $A$ is symmetric, it is diagonalizable by an orthogonal transformation, and we may assume that $A = \text{diag}(\lambda_1, \dots, \lambda_n)$.
With the notation $\lambda = \lambda_j$, it follows that $|\be^T\!p(A)\be| = |\sum_{i\ne j} e_i^2\, p(\lambda_i)|$, from which the result follows easily. 
\end{proof}
The above expressions imply the following lower and upper bounds on the sensitivity of the approximate eigenvalues as function of the approximate eigenvector.
The assumption $\lambda \ne 0$ in the next propositions may be viewed as non-restrictive: for zero eigenvalues we can consider the absolute error $|\theta(\bu)-\lambda|$ instead.

\medskip
\begin{proposition} \label{prop:1}
Suppose $\bu = \bx + \be$ is an approximate eigenvector corresponding to a simple eigenvalue $\lambda \ne 0$ of a symmetric $A$, with $\|\bx\|=1$, $\be \perp \bx$, and $\eps = \|\be\|$. Then, up to $\calo(\eps^4)$-terms, for the sensitivity of the Rayleigh quotient (as function of $\bu$) it holds that
\begin{equation} \label{eq:brq}
\min_{\lambda_i \ne \lambda} \frac{|\lambda_i-\lambda|}{|\lambda|} \ \eps^2 \ \lesssim \
\frac{|\theta(\bu)-\lambda|}{|\lambda|}
\ \lesssim \ \max_{\lambda_i \ne \lambda} \frac{|\lambda_i-\lambda|}{|\lambda|} \ \eps^2.
\end{equation}
If $\tau$ is not equal to an eigenvalue, then for the harmonic Rayleigh quotient with target $\wt \theta_\tau$:
\begin{equation}
\label{eq:bhrq}
\min_{\lambda_i \ne \lambda} \frac{\vert(\lambda_i-\tau)(\lambda_i-\lambda)\vert}{\vert\lambda(\lambda-\tau)\vert}\, \eps^2 \ \lesssim \
\frac{|\wt \theta_\tau(\bu)-\lambda|}{|\lambda|}
\ \lesssim \ \max_{\lambda_i \ne \lambda} \frac{\vert(\lambda_i-\tau)(\lambda_i-\lambda)\vert}{\vert\lambda(\lambda-\tau)\vert}\, \eps^2.
\end{equation}
In particular, for the sensitivity of the harmonic Rayleigh quotient (i.e., $\tau=0$) we have 
\[
\min_{\lambda_i \ne \lambda} \frac{\vert\lambda_i(\lambda_i-\lambda)\vert}{\lambda^2}\, \eps^2\ \lesssim \
\frac{|\wt \theta(\bu)-\lambda|}{|\lambda|}
\ \lesssim \ \max_{\lambda_i \ne \lambda} \frac{\vert\lambda_i(\lambda_i-\lambda)\vert}{\lambda^2}\, \eps^2.
\]
\end{proposition}
\begin{proof}
This follows from \eqref{rq-sens}, \eqref{hrq-sens}, and \eqref{hrq-tau-sens}, using Lemma~\ref{lemma} with $p(t) = t-\lambda$ for \eqref{rq-sens} and $p(t) = t\,(t-\lambda)$ for \eqref{hrq-tau-sens}. The bound for \eqref{rq-sens} is a specific case of \eqref{hrq-sens}, when $\tau = 0$.
\end{proof}

We mention that the bound \eqref{eq:brq} has a slightly improved version, as follows (cf., e.g., \cite[Thm.~2.1]{sleijpen2003optimal} for the smallest eigenvalue).
For this context we introduce $\wt\be = \frac{\be}{\|\be\|}$ of unit length, $\bu = \bx + \be$ as before, and $\wt\bu = \frac{\bu}{\|\bu\|}$, where $\|\bu\|^2 = 1+\eps^2$.
We can decompose $\wt\bu = \cos(\wt\bu,\bx) \, \bx + \sin(\wt\bu,\bx) \, \wt\be$. Then an easy computation gives $\wt\bu^T\!A\wt\bu-\lambda = \sin^2(\wt\bu,\bx) \ \wt\be^T\!(A-\lambda I)\,\wt\be$ and therefore 
\begin{equation} \label{sharp-sin2}
|\wt\bu^T\!A\wt\bu-\lambda| \le \max_{\lambda_i \ne \lambda} |\lambda_i-\lambda| \, \sin^2(\wt\bu,\bx).
\end{equation}
To connect the approximate bound \eqref{eq:brq} to \eqref{sharp-sin2}, we note that $\sin^2(\wt\bu,\bx) = \frac{\eps^2}{1+\eps^2}$, which is asymptotically equal to $\eps^2$. The bound \eqref{sharp-sin2} is exact (i.e., not an asymptotic bound as in Proposition~\ref{prop:1}) and sharp, but it is asymptotically equal to our expression.

An exact bound for the error of the harmonic Rayleigh quotient with target is provided in \cite[Thm.~5.2]{sleijpen2003use}. Under certain assumptions mentioned in \cite{sleijpen2003use} for the target $\tau$ and $\eps$, with $\bx = \bx + \eps\,\wt\be$, it is shown that
\[
\frac{\bu^TA(A-\tau I)\,\bu}{\bu^T(A-\tau I)\,\bu} - \lambda = \eps^2\frac{\wt\be^T(A-\lambda I)(A - \tau I)\,\wt\be}{(\lambda - \tau) + \eps^2 \, (\wt\be^TA\wt\be - \tau)} \le \eps^2 \, \max_{\lambda_i \ne \lambda}\frac{(\lambda_i-\lambda)(\lambda_i - \tau)}{(\lambda - \tau) + \eps^2(\lambda_i - \tau)}.
\]
Since we do not make any assumptions on the target or the value of $\|\be\|$, our upper bound \eqref{eq:bhrq} includes absolute values. Discarding the $\eps^2$-term in the denominator yields an approximate upper bound accurate to $\calo(\eps^4)$-terms. In conclusion, \eqref{eq:bhrq} is less sharp in some situations but more general, and asymptotically the same as the result of \cite{sleijpen2003use}. Another good reason to consider \eqref{eq:brq} and \eqref{eq:bhrq} is that approximate bounds for the error of the homogeneous Rayleigh quotient can also be obtained via the second-order approximation of the quotient, as we have done in Proposition~\ref{prop:1}. At the end of this section, we will also point out that the three approximate bounds have the same form.

\subsection{Homogeneous Rayleigh quotient} 
\label{sec:homobounds}
We now determine approximate bounds for the sensitivity of the homogeneous Rayleigh quotient.
\begin{proposition}
Suppose $\bu = \bx + \be$ is an approximate eigenvector corresponding to a simple eigenvalue $\lambda \ne 0$ of a symmetric $A$, with $\|\bx\|=1$, $\be \perp \bx$, and $\eps = \|\be\|$. Then, up to $\calo(\eps^4)$-terms, for the sensitivity of the homogeneous Rayleigh quotient (as function of $\bu$) it holds that:
\[
\frac{\eps^2}{\lambda^2+1} \, \min_{\lambda_i \ne \lambda} \vert(\lambda_i-\lambda)(\lambda_i+\lambda^{-1})\vert \, \lesssim \, \frac{|\alpha(\bu)-\lambda|}{|\lambda|}
\, \lesssim \, \frac{\eps^2}{\lambda^2+1} \, \max_{\lambda_i \ne \lambda} \vert(\lambda_i-\lambda)(\lambda_i+\lambda^{-1})\vert.
\]
\end{proposition}
\begin{proof}
From Proposition~\ref{prop:prop}(iii), the homogeneous Rayleigh quotient can be written as
\[
\alpha = \sqrt{\aux^2+1} - \aux, \qquad \aux = (2\ \bu^T\!A\bu)^{-1} \, (\bu^T\bu-\bu^T\!A^2\bu).
\]
We can express $\aux$ in terms of $\lambda$ and $\be$ as
\begin{equation*}
\aux = \frac{1-\lambda^2 + \|\be\|^2 - \|A\be\|^2}{2\,(\lambda + \be^TA\be)}.
\end{equation*}
Furthermore, when $\|\be\| = \eps \to 0$, 
\[
\aux = \tfrac12 \, \lambda^{-1} \, (1 - \lambda^{-1} \, \be^T\!A\be)\,(1 - \lambda^2 + \|\be\|^2 -\|A\be\|^2 )+ \calo(\eps^4).
\]
For the square root, we use the first-order approximation $\sqrt{1+t} = 1 + \frac12 t + \calo(t^2)$. Discarding $\calo(\eps^4)$-terms yields
\begin{align*}
\sqrt{\aux^2+1} &= \tfrac12 \, \lambda^{-1} \, (1 + \lambda^{-1} \, \be^T\!A\be)^{-1} \sqrt{(1-\lambda^2 + \|\be\|^2 - \|A\be\|^2)^2 + 4\,(\lambda + \be^TA\be)^2} \\
&= \tfrac12 \lambda^{-1} (1 - \lambda^{-1} \be^TA\be) \big(1 + \lambda^2 + \tfrac{1-\lambda^2}{1+\lambda^2}(\|\be\|^2 -\|A\be\|^2) + \tfrac{4\lambda}{1+\lambda^2}\,\be^T\!A\be \big).
\end{align*}
Simple computations show that
\begin{align}
\label{eq:homoapprox}
\alpha & = (1 - \lambda^{-1} \be^TA\be) \Big(\lambda -\frac{\lambda}{\lambda^2+1}(\|\be\|^2 -\|A\be\|^2) +\frac{2}{\lambda^2+1}\be^T\!A\be \Big) + \calo(\eps^4)\nonumber\\
& = \lambda +\frac{\lambda}{\lambda^2+1} \, \be^T(A + \lambda^{-1}I)(A-\lambda I) \, \be + \calo(\eps^4).
\end{align}
The thesis follows from Lemma~\ref{lemma}, with the polynomial $p(t) = (t + \lambda^{-1})(t - \lambda)$.

An alternative method to derive this result is via the Implicit Function Theorem. From \eqref{eq:homogal}, for the pair $(\alpha(\bu), \bu)$, we have the implicit equation $F(\alpha(\bu), \bu) := \bu^T (A-\alpha I)(A+\alpha^{-1}I) \, \bu = 0$.
Then we know that
\[
\nabla \alpha(\bu) = -\Big( \frac{\partial F}{\partial \alpha} \Big)^{-1} \, \frac{\partial F}{\partial \bu} \\
= 2 \ [\bu^T(A+\alpha^{-2}A)\,\bu]^{-1} \ (A-\alpha I)(A+\alpha^{-1} I)\, \bu.
\]
Since $\nabla\alpha(\bx) = \zero$, we consider the second-order approximation of the homogeneous Rayleigh quotient $\alpha(\bx+\be) \approx \alpha(\bx) + \frac12 \, \be^T \nabla^2 \alpha(\bx) \, \be$.
The action of $\nabla^2 \alpha(\bx)$ can be obtained from the first-order approximation of the gradient in $\bx$:
\begin{align*}
\nabla^2\alpha(\bx)\,\be & = \nabla \alpha(\bx+\be) - \nabla\alpha(\bx) + \calo(\eps^2) \\
& = \frac{2 \, \lambda}{1+\lambda^2} \ (A-\lambda I)(A+\lambda^{-1} I)\,\be + \calo(\eps^2).
\end{align*}
This also leads to the expression in \eqref{eq:homoapprox}.
\end{proof}

Interestingly, by replacing $\tau = -\lambda^{-1}$ in the approximation of the harmonic Rayleigh quotient with target \eqref{hrq-tau-sens}, we get the approximation \eqref{eq:homoapprox} for the homogeneous Rayleigh quotient.

We now compare the different asymptotic bounds. With a little abuse of notation, let $\theta(\bu)$ be any of the three Rayleigh quotients in $\bu$. Then the previous results can be summarized by 
\begin{equation*}
\frac{\eps^2}{\vert\lambda\vert}\ \min_{\lambda_i \ne \lambda} |p_\lambda(\lambda_i)| \ \ \lesssim \ \
\frac{|\theta(\bu)-\lambda|}{|\lambda|}
\ \ \lesssim \ \ \frac{\eps^2}{\vert\lambda\vert}\ \max_{\lambda_i \ne \lambda} |p_\lambda(\lambda_i)|, 
\end{equation*}
where the different polynomials $p_\lambda(\cdot)$ are reported in Table~\ref{tab:poly}, along with their pointwise asymptotic behavior for $\lambda\to 0$ and $\lambda\to\pm\infty$. This table provides a naive explanation of what we will observe in the experiments: the homogeneous Rayleigh quotient tends to have a similar behavior to the Rayleigh quotient for small eigenvalues, and it follows the harmonic Rayleigh quotient for large eigenvalues.
\begin{table}[ht!]
\footnotesize
\centering
\caption{Polynomials $p_{\lambda}(\cdot)$ appearing in the asymptotic bounds for the relative error of the various Rayleigh quotients.}
\begin{tabular}{lccc}
\toprule
& $p_\lambda(t)$ & $\lambda\to 0$ & $\lambda\to\pm\infty$\\
\midrule
RQ & $t-\lambda$ & $\sim t$ & $\sim -\lambda$ \\[1mm]
Harmonic RQ & $\frac{1}{\lambda} \, t \, (t-\lambda)$& $\sim \frac{t^2}{\lambda}$ & $\sim -t$\\[1mm]
Homogeneous RQ & $\frac{1}{\lambda^2+1} (t-\lambda)(\lambda t + 1)$ &  $\sim t$ & $\sim -t$\\
\bottomrule
\end{tabular}
\label{tab:poly}
\end{table}

We remark that, for the homogeneous Rayleigh quotient, $p_{\lambda} = p_{-\lambda^{-1}}$.
In the special case that both $\lambda$ and its anti-reciprocal $-\lambda^{-1}$ lie in the spectrum of $A$, the asymptotic lower bound for the homogeneous Rayleigh quotient will be zero, since in general $\lambda \ne -\lambda^{-1}$ and both eigenvalues are roots of $p_\lambda$.

For the Rayleigh quotient, we have that $\max_{\lambda_i \ne \lambda}\vert\lambda_i - \lambda\vert \in \{\vert\lambda - \lambda_1\vert, \, \vert\lambda - \lambda_n\vert\}$. In contrast, the upper bounds for the harmonic Rayleigh quotient or the homogeneous Rayleigh quotient are more complicated. Viewing $p_\lambda$ as a continuous function on $[\lambda_1, \lambda_n]$, the maximum of the corresponding two upper bounds may be attained at the vertex of the parabola $p_\lambda$, or at the boundary points $\lambda_1$ and $\lambda_n$. Thus, in view of the discrete spectrum, we know that
\begin{equation}
\label{eq:maxmax}
\max_{\lambda_i \ne \lambda} |p_\lambda(\lambda_i)| \le \max_{\lambda_i \in \{\lambda_1, \ \lambda_n, \ \overline \lambda\}} |p_\lambda(\lambda_i)|,
\end{equation}
where $\overline\lambda = \frac12 \lambda$ for the harmonic Rayleigh quotient, while $\overline \lambda = \frac12 (\lambda - \lambda^{-1})$ for the homogeneous Rayleigh quotient. Although $p_\lambda(t)$ always has a factor of $t-\lambda$, due to the quadratic nature of the upper bound for the harmonic and the homogeneous Rayleigh quotient, we cannot make any a priori comparison with the bound for the Rayleigh quotient; without additional information about the spectrum, we cannot improve on these considerations.
The behavior of the relative errors and upper bounds for the various Rayleigh quotients is shown in the next section.

\subsection{Comparison of the Rayleigh quotients}
\label{sec:sensexp}
Without loss of generality, we may assume that $A = \text{diag}(\lambda_1,\dots,\lambda_n)$. We consider a family of $100 \times 100$ positive definite diagonal matrices where the eigenvalues have uniform distribution on the interval $[0,2\sigma]$, i.e., $\lambda_i \sim U(0, 2\sigma)$, $\sigma > 0$. The second family contains $100 \times 100$ indefinite diagonal matrices where the eigenvalues have Gaussian distribution with mean $0$ and variance $\sigma^2$, which we indicate by $\lambda_i \sim \mathcal{N}(0, \sigma)$, $\sigma > 0$. Since the eigenvalues are drawn from continuous probability distributions, there is zero probability that two eigenvalues are equal, so the eigenvectors are well defined with probability one.

An eigenpair of a generic diagonal $A$ is $(\lambda, \bx)$, where $\bx$ is one of the vectors of the canonical basis of $\mathbb{R}^{100}$.
To compute a random perturbation $\bu = \bx + \eps\,\be$, such that $\|\be\| = 1$ and $\bx \perp \be$, we start from a vector with random Gaussian components, project it onto the orthogonal complement of $\bx$ and normalize to get $\be$. Finally we normalize $\bu$.

\begin{figure}[ht]
\centering
\includegraphics[width=\textwidth]{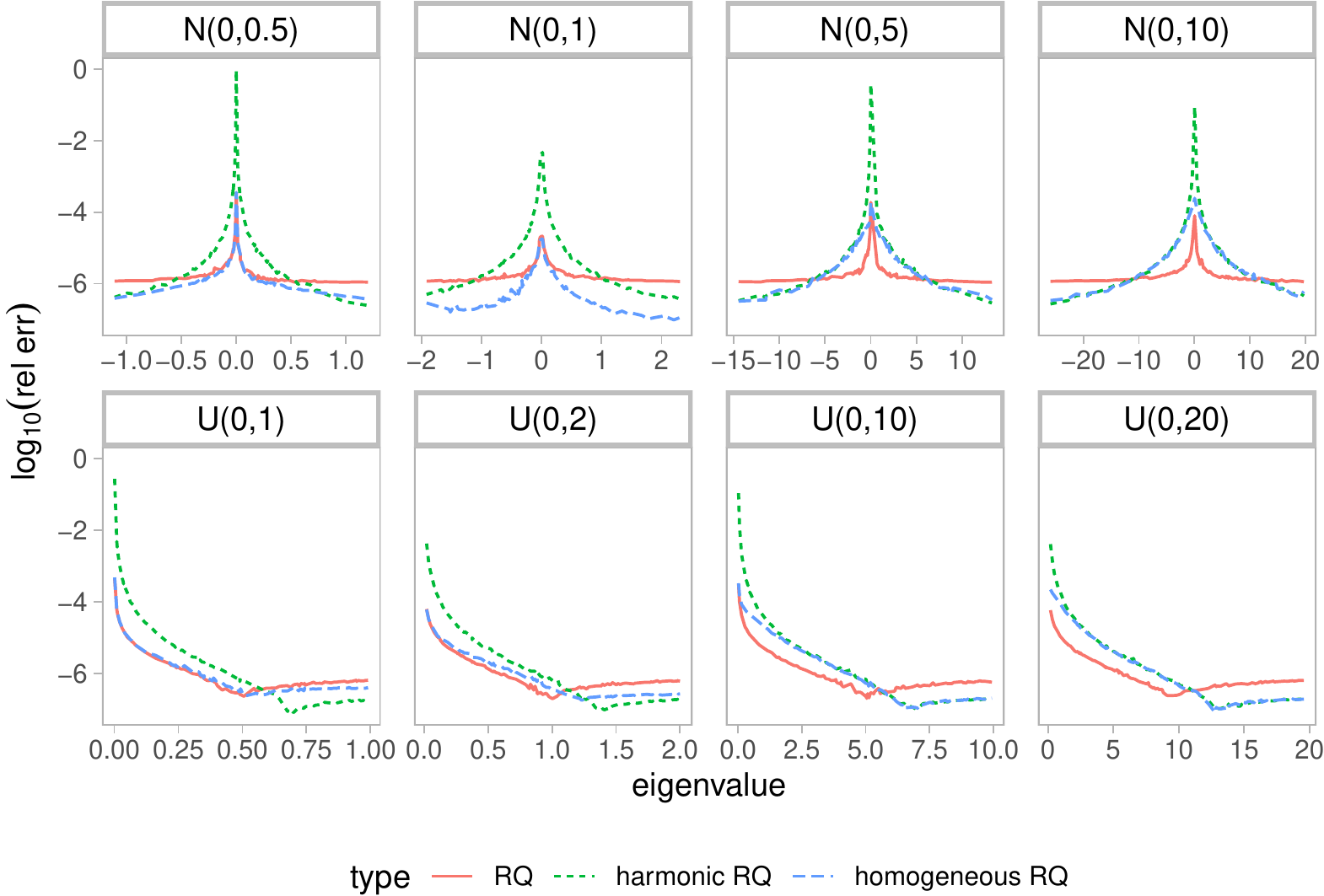}
\caption{Sensitivity of Rayleigh quotients (RQ) for different spectra. The label of each graph indicates the distribution of the $100$ eigenvalues. Each eigenvalue is plotted against its corresponding maximum relative error, in logarithmic scale.}
\label{fig:sigmas}
\end{figure}

First, we study the behavior of the different Rayleigh quotients when perturbed, with $\sigma \in \{0.5, 1, 5, 10\}$ and $\eps = 0.001$. For each eigenvalue $\bx$ we draw $100$ random perturbations $\bu_{j}$, and take the maximum relative error. The quantity of interest for, e.g., the standard Rayleigh quotient is
\[
\max_j \frac{\vert\theta(\bu_{j})-\lambda\vert}{\vert\lambda\vert}.
\]
Figure~\ref{fig:sigmas} shows the maximum relative error as a function of the estimated eigenvalue $\lambda$.
In these examples, there are peaks close to zero since we consider relative errors; therefore, the Rayleigh quotients are more sensitive when small eigenvalues are approximated. In particular, the harmonic Rayleigh quotient is always more sensitive than the other two for small eigenvalues. By looking at the two families separately, Gaussian and uniformly distributed, we notice that as $\sigma$ increases, the curves of the Rayleigh quotient and the harmonic Rayleigh quotient remain almost unchanged.
In the uniformly distributed case, when $\sigma = 0.5$, the homogeneous Rayleigh quotient closely follows the Rayleigh quotient in the first half of the spectrum and slightly departs in the second half. Starting from $\sigma = 1$, the curve gets increasingly closer to the harmonic Rayleigh quotient until the two curves are very similar, except for the smaller eigenvalues, where the relative error of the homogeneous Rayleigh quotient is lower than the one of the harmonic Rayleigh quotient.
This trend can be partially explained by the asymptotic behavior of the polynomials in Table~\ref{tab:poly}: the homogeneous and \gf{the} harmonic Rayleigh quotient\gf{s} have the same behavior for large eigenvalues, while the homogeneous Rayleigh quotient tends to be similar to the standard Rayleigh quotient for small eigenvalues.

We observe the same characteristics in the Gaussian family, with the interesting addition that when $\sigma = 1$, the homogeneous Rayleigh quotient is even less sensitive than the harmonic Rayleigh quotient.
We have also considered different magnitudes of perturbation $\eps$, where we  observe that the relative errors increase as $\eps$ increases, but the relative positions of the Rayleigh quotients remain the same.

\begin{figure}[ht] 
\centering
\subfloat[Eigenvalues drawn from $\mathcal{N}(0,1)$]{
\includegraphics[width = 0.85\textwidth]{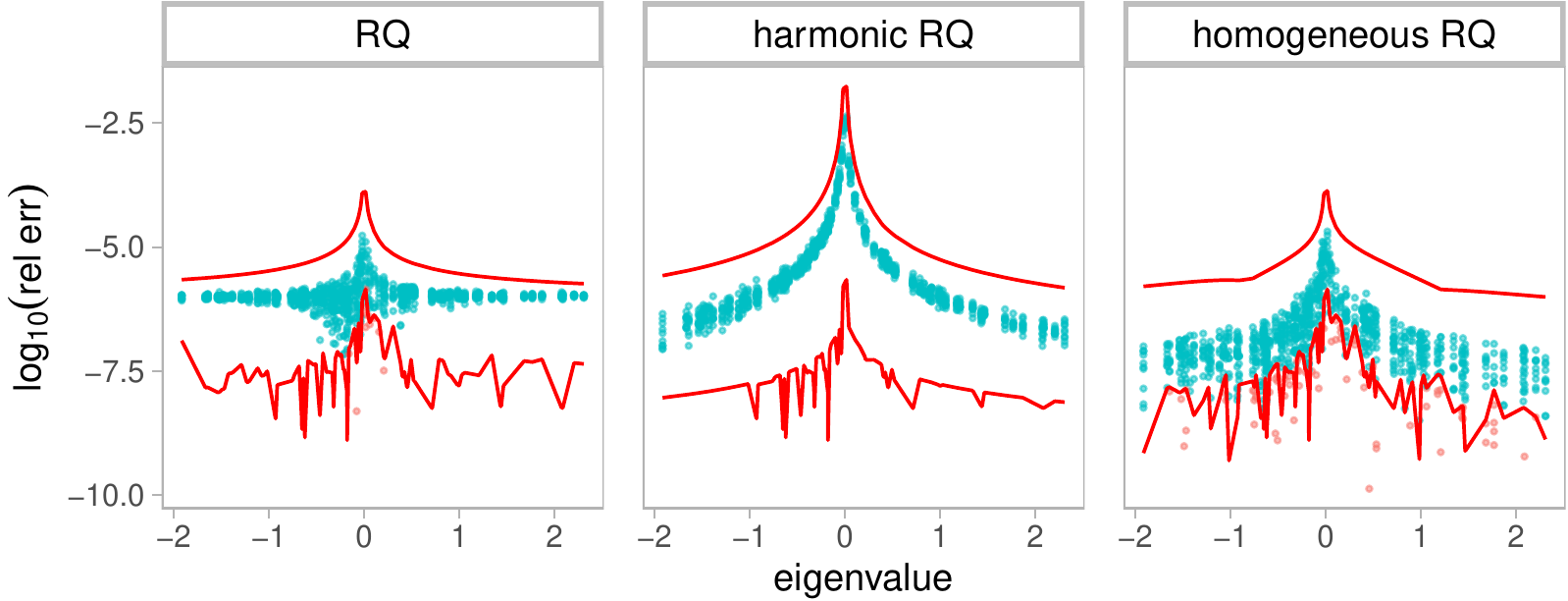}}

\subfloat[Eigenvalues drawn from $U(0,2)$]{
\includegraphics[width = 0.85\textwidth]{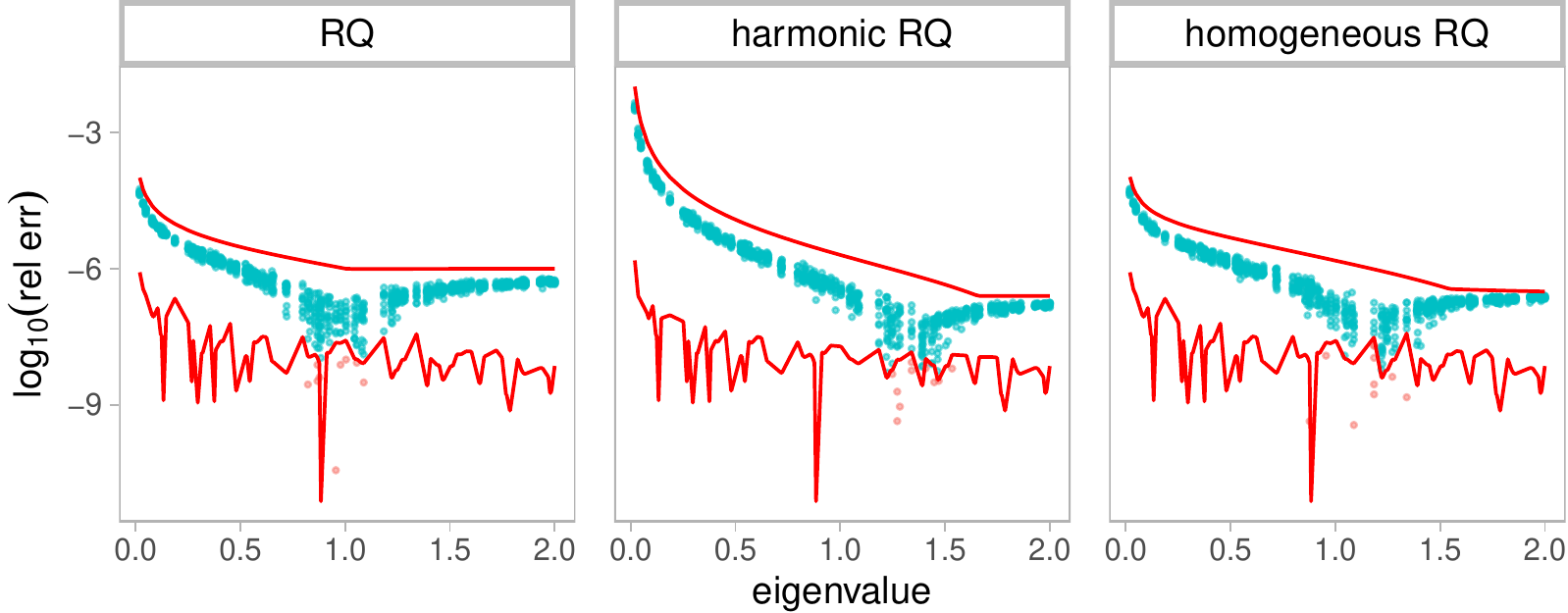}}
\caption{Sensitivity of the Rayleigh quotients for two different spectra. Each eigenvalue is plotted against their $10$ corresponding relative errors, one for each perturbation of size $\eps = 0.001$. The asymptotic upper and lower bounds are shown as red lines. Points that fall outside these bounds are plotted in red.}
\label{fig:norm-and-unif}
\end{figure}

Finally, we look at the asymptotic bounds of Section~\ref{sec:sens}, when $\lambda_i \sim \mathcal{N}(0,1)$ and $\lambda_i \sim U(0,2)$. Instead of just showing the maximum relative error, we draw $10$ random perturbations for each eigenvector, and plot the corresponding relative errors. The results are presented in Figure~\ref{fig:norm-and-unif}. 
We observe that the lower bounds are more erratic than the upper bounds. This can be explained as follows. In Section~\ref{sec:homobounds} we have remarked that all bounds depend on a continuous nonnegative function $\vert p_\lambda(t)\vert$, which has a zero in $t = \lambda$ (cf.~Table~\ref{tab:poly}). Therefore, given a $\lambda_j$, the minimizer $\argmin_{\lambda_i \ne \lambda} \vert p_\lambda(\lambda_i)\vert$ is typically either $\lambda_{j-1}$ or $\lambda_{j+1}$; the maximizer is generally either $\lambda_1$ or $\lambda_n$ (cf.~Table~\ref{tab:poly}). For this reason, the lower bound shows more variability. As already mentioned, while this statement holds for the Rayleigh quotient, this is not always true for the harmonic and the homogeneous Rayleigh quotients. 

Regarding the behavior of the relative errors, the plots in Figure~\ref{fig:norm-and-unif} reasonably reflect the ones in Figure~\ref{fig:sigmas}. For the uniformly distributed eigenvalues, the upper bound seems to be sharp at the extremes of the spectrum, while for the Gaussian family it is sharper only close to the smallest eigenvalues (in magnitude). The upper bound for the Rayleigh quotient is also sharp for the largest eigenvalues. While the lower bound of the uniformly distributed family poorly reflects the behavior of the relative error, in the Gaussian family it is tighter. In particular, it is very accurate for the homogeneous Rayleigh quotient. In Figure~\ref{fig:sigmas}, we have already remarked that the maximum relative error in the homogeneous Rayleigh quotient is lower than for the other Rayleigh quotients. Now we see that for some perturbations, the sensitivity can be even much lower.

So far we have discussed the properties of the homogeneous Rayleigh quotient and compared it to the well-known standard and harmonic Rayleigh quotient. Now we extend the homogeneous Rayleigh quotient for the generalized eigenvalue problem and present an application of Rayleigh quotients to gradient methods for unconstrained optimization problems.

\section{The generalized eigenvalue problem} \label{sec:gep}
Let us consider the generalized eigenvalue problem (GEP) $A\bx = \lambda B\,\bx$, where $A$ is symmetric (definite or indefinite), and $B$ is SPD. This problem has $n$ real eigenvalues $\lambda_1 \le \cdots \le \lambda_n$.
Although a common generalized Rayleigh quotient is $\frac{\bu^T\!A\bu}{\bu^T\!B\,\bu}$ (see, e.g., \cite[Ch.~15]{Par98}), it turns out that, in our context, the most relevant quantities are
\begin{equation} \label{rqgep}
\theta = \frac{\bu^T\!AB\,\bu}{\bu^T\!B^2\bu}, \qquad
\wt \theta = \frac{\bu^T\!A^2\bu}{\bu^T\!AB\,\bu},
\end{equation}
for nonzero $\bu^T\!AB\,\bu$.
The Rayleigh quotient $\theta$ satisfies $A\bu-\theta B\,\bu \perp B\,\bu$, and, equivalently, it is the solution to $\min_{\gamma} \|A\bu-\gamma B\,\bu\|$. 
The harmonic Rayleigh quotient $\wt \theta$ satisfies $A\bu-\theta B\,\bu \perp A\bu$ and solves $\min_{\gamma} \|\gamma^{-1} A\bu-B\,\bu\|$.
The extension of the homogeneous Rayleigh quotient for the GEP is rather straightforward.
It is the solution to
\begin{equation} \label{gen-homo-sec}
\min_{(\alpha_1,\,\alpha_2) \, \in \, \mathbb P} \ 
\frac{\|\alpha_1 B\, \bu - \alpha_2 \, A\bu\|}{\gf{\sqrt{\alpha_1^2+\alpha_2^2}}}
= \min_{\alpha_1^2+\alpha_2^2=1} \ \|[B\,\bu \ \ - \!\!A\bu] \smtxa{c}{\alpha_1 \\ \alpha_2}\|.
\end{equation}
As for the standard eigenvalue problem, this amounts to solve a reduced SVD of an $n \times 2$ matrix, or an eigenvalue problem involving a $2 \times 2$ matrix.

The next result is a generalization of Proposition~\ref{prop:prop}.
Let $C = [B\,\bu \ \ -\!\!A\bu]$. 

\begin{proposition} 
Suppose that $\bu^T\!AB\,\bu \ne 0$ and denote by $\mu$ the smallest eigenvalue of $C^T\!C$. Then the following properties hold.
\begin{itemize}
\item[(i)] $\mu$ is a simple eigenvalue of $C^T\!C$, or equivalently $\sqrt{\mu}$ is a simple singular value of $C$; its corresponding eigenvector $[\alpha_1, \, \alpha_2]^T$ is the unique minimizer of \eqref{gen-homo-sec}, up to orthogonal transformations. 
\item[(ii)] $\mu = \tfrac12 \, \big[\bu^T\!A^2\bu+\bu^T\!B^2\bu - \sqrt{\smash[b]{(\bu^T\!A^2\bu-\bu^T\!B^2\bu)^2+4\,(\bu^T\!AB\,\bu)^2}} \ \big].$
\item[(iii)] For the generalized homogeneous Rayleigh quotient $\alpha$ we have
\begin{align}
\alpha &= 
 (2\ \bu^T\!AB\,\bu)^{-1} \, \big[ \,\bu^T\!A^2\,\bu - \bu^T\!B^2\bu + \sqrt{(\bu^T\!A^2\,\bu - \bu^T\!B^2\bu)^2 + 4\,(\bu^T\!AB\,\bu)^2} \,\big].
\end{align}
\item[(iv)] $\mu=0$ if and only if $\bu$ is an eigenvector. If $\bu$ is an eigenvector, then the corresponding homogeneous Rayleigh quotient $\alpha$ is the corresponding eigenvalue.
\item[(v)] $0 \le \mu < \min(\bu^T\!A^2\bu, \ \bu^T\!B^2\bu)$.
\item[(vi)] We have the following inequalities
\begin{equation*}
\begin{cases}
\theta \le \alpha \le \wt \theta &\text{if}\quad\bu^T\!AB\,\bu > 0, \\
\wt \theta \le \alpha \le \theta &\text{if}\quad\bu^T\!AB\,\bu < 0.
\end{cases}
\end{equation*}
\end{itemize}
\end{proposition}
\begin{proof}
This proof follows the exact same lines as those of Proposition~\ref{prop:prop}.
\end{proof}
As in Section~\ref{sec:homo}, we can show that the generalized homogeneous Rayleigh quotient is a solution to a quadratic equation.
\begin{proposition}\label{prop:genquadeq}
Let $\bu^T\!AB\,\bu \ne 0$. The generalized homogeneous Rayleigh quotient \eqref{eq:homorq} is the solution to 
\begin{equation}
\label{eq:genhomopoly}
(\bu^T\!AB\,\bu) \, \alpha^2 + (\bu^T\!B^2\bu-\bu^T\!A^2\bu) \, \alpha - (\bu^T\!AB\,\bu) = 0,
\end{equation}
which satisfies $\alpha\ \bu^T\!AB\,\bu > 0$.
\end{proposition}
\begin{proof}
The proof is similar to that of Proposition~\ref{prop:quadeq}. 
\end{proof}

We derive a Galerkin condition from Proposition~\ref{prop:genquadeq} as follows. Since $\bu^T\!AB\,\bu = \bu^T\!B\,A\bu$, \eqref{eq:genhomopoly} is equivalent to $\bu^T\!(\alpha A + B)(A - \alpha B)\,\bu = 0$.
Therefore, we can write the corresponding {\em nonlinear} Galerkin condition in homogeneous coordinates as
\begin{equation*}
(\alpha_2A - \alpha_1B)\,\bu \perp (\alpha_1A + \alpha_2B)\,\bu.
\end{equation*}
Finally, we discuss a bound for the chordal metric, analogous to \eqref{eq:chord}. For an SPD $B$, as a generalization to \cite[Thm.~4.5.1]{Par98}, one can show that there exists a generalized eigenvalue $\lambda = \lambda(A,B)$ of $(A, B)$ such that $\vert\lambda - \alpha\vert \le (\lambda_1(B)\|\bu\|)^{-1}\|(A-\alpha\, B)\,\bu\|$, where $\lambda_1(B)$ is the smallest eigenvalue of $B$.
Indeed, if $\alpha$ is not an eigenvalue of $(A,B)$,
\begin{align*}
\|\bu\| &= \|(A-\alpha\, B)^{-1}\,(A-\alpha\, B)\,\bu\| \ge \lambda_1(B)\,\|(B^{-1}A-\alpha\, I)^{-1}\| \, \|(A-\alpha\, B)\,\bu\|\\
&\ge \lambda_1(B) \, \min_i\vert \lambda_i(A,B) - \alpha\vert\cdot\|(A-\alpha\, B)\,\bu\|.
\end{align*}
Therefore, for the chordal distance between $\lambda$ and $\alpha$ we have
\[
\frac{\vert\lambda - \alpha\vert}{\sqrt{1+\lambda^2} \, \sqrt{1+\alpha^2}} \le \frac{1}{\lambda_1(B)}\cdot\frac{1}{\sqrt{1+\lambda^2} \, \|\bu\|}\cdot\frac{\|(A-\alpha\, B)\,\bu\|}{\sqrt{1+\alpha^2}}.
\]
This upper bound is minimized for $\alpha$ equal to the generalized homogeneous Rayleigh quotient. We notice that this inequality differs from \eqref{eq:chord} by the factor $(\lambda_1(B))^{-1}$. 

In the next section, we return to the standard homogeneous Rayleigh quotient, and study its use as stepsize in gradient methods.

\section{A homogeneous stepsize for gradient methods} \label{sec:step}
In gradient methods for nonlinear optimization, inverse Rayleigh quotients are popular choices for stepsizes. We refer to \cite{daniela2018steplength} for a nice recent review about steplength selection.
We exploit the inverse of the homogeneous Rayleigh quotient as a new \emph{homogeneous stepsize} for gradient-type methods.
Consider the unconstrained minimization of a smooth function $f: \R^n \to \R$,
\begin{equation*}
\min_{\bx \in \mathbb{R}^n} \, f(\bx).
\end{equation*}
Gradient methods are of the form
\begin{equation*} 
\bx_{k+1} = \bx_k - \beta_k \, \bg_k = \bx_k - \alpha_k^{-1} \, \bg_k,
\end{equation*}
where $\bg_k = \nabla f(\bx_k)$, $\beta_k > 0$ is the stepsize, and $\alpha_k$ the inverse stepsize.
As usual we write $\bs_{k-1} = \bx_k-\bx_{k-1}$ and $\by_{k-1} = \bg_k-\bg_{k-1}$. 

In \cite{bb1988}, the Barzilai--Borwein stepsizes
\begin{equation} \label{bb1}
\beta_k^{\rm{BB1}} = \frac{\bs_{k-1}^T \, \bs_{k-1}}{\bs_{k-1}^T \, \by_{k-1}}, \qquad
\beta_k^{\rm{BB2}} = \frac{\bs_{k-1}^T \, \by_{k-1}}{\by_{k-1}^T \, \by_{k-1}}
\end{equation}
have been introduced. The motivation is that, when we approximate the Hessian in $\bx_{k-1}$ by a scalar multiple of the identity, i.e., $\nabla^2 f(\bx_{k-1}) \approx \gamma\, I$, the corresponding inverse stepsizes satisfy the following least squares secant conditions:
\begin{equation} \label{bbsec}
\alpha_k^{\rm{BB1}} = \argmin_{\gamma} \|\by_{k-1} - \gamma\,\bs_{k-1}\|, \qquad
\alpha_k^{\rm{BB2}} = \argmin_{\gamma} \|\gamma^{-1} \, \by_{k-1}-\bs_{k-1}\|,
\end{equation}
with $\alpha_k^{\rm{BB1}} = (\beta_k^{\rm{BB1}})^{-1}$ and $\alpha_k^{\rm{BB2}} = (\beta_k^{\rm{BB2}})^{-1}$. Moreover, both stepsizes \eqref{bb1} can be seen as inverse Rayleigh quotients of a certain matrix $H_k$ at $\bs_{k-1}$. In fact, for any $f$ it holds
\begin{equation*}
\by_{k-1} = \nabla f(\bx_k) - \nabla f(\bx_{k-1}) = H_k \, (\bx_k-\bx_{k-1}) = H_k \, \bs_{k-1},
\end{equation*}
where $H_k := \int_0^1 \nabla^2 f((1-t) \, \bx_{k-1}+ t \, \bx_k) \, dt$ is an average Hessian on the line piece between $\bx_{k-1}$ and $\bx_k$. 
From this relation, it is easy to see that the minimum residual conditions \eqref{opt1b} and \eqref{opt2b} are equivalent to the secant conditions \eqref{bbsec} for $\bu = \bs_{k-1}$ and $A = H_k$.

We introduce the homogeneous BB stepsize (HBB) as the inverse of the homogeneous Rayleigh quotient of $H_k$ in $\bs_{k-1}$. The HBB stepsize is given by the quotient $\beta_k^{\rm HBB} = \alpha_{2,k} / \alpha_{1,k}$, where the pair $(\alpha_{1,k}$, $\alpha_{2,k})$ solves
\[
\argmin_{\alpha_1^2+\alpha^2_2 = 1} \ \|\alpha_1 \, \bs_{k-1}-\alpha_2 \, \by_{k-1}\|.
\]
As in \eqref{bbsec}, this is equivalent to the minimum residual condition \eqref{homo-sec} for $\bu = \bs_{k-1}$ and $A = H_k$. 
Therefore (cf. Proposition~\ref{prop:prop})
\begin{equation}
\label{eq:homobb}
\beta_k^{\rm HBB} = \frac{\|\bs_{k-1}\|^2-\|\by_{k-1}\|^2 + \sqrt{\smash[b]{(\|\bs_{k-1}\|^2-\|\by_{k-1}\|^2)^2 + 4\,(\bs_{k-1}^T\by_{k-1})^2}}}{2\ \bs_{k-1}^T\by_{k-1}}.
\end{equation}
\gf{As mentioned in Section~\ref{sec:intro}, this stepsize has been proposed first by Li, Zhang, and Xia \cite{LZX22,LX21}, obtained from a total least squares secant condition.}
In this section, we will also introduce an alternating variant of this HBB step.

In Section~\ref{sec:exp} we will carry out some experiments to test the behavior of the HBB stepsize when plugged into a gradient method for general differentiable functions. A pseudocode of this method is provided in Algorithm~\ref{algo:nonlin}. 
\begin{algorithm}
\caption{A homogeneous gradient method for minimization of general functions}
\label{algo:nonlin}
{\bf Input}: Continuous differentiable function $f$, initial guess $\bx_0$, initial stepsize $\beta_0 > 0$, tolerance {\sf tol}; safeguarding parameters $\beta_{\max} > \beta_{\min} > 0$; line search parameters $c_{\rm{ls}}$, $\sigma_{\rm{ls}} \in (0,1)$; memory integer $M>0$\\
{\bf Output}: Approximation to minimizer $\argmin_{\bx} f(\bx)$ \\
\begin{tabular}{rl}
{\footnotesize 1}: & Set $\bg_0 = \nabla f(\bx_0)$ \\
& {\bf for} $k = 0, 1, \dots$ \\
{\footnotesize 2}: & \phantom{M} $\nu_k = \beta_k$, \ $f_{\text{ref}} = \max \, \{ \, f(\bx_{k-j}) \, : \, 0 \le j \le \min(k,M-1) \, \}$ \\
{\footnotesize 3}: & \phantom{M} {\bf while} \ $f(\bx_k-\nu_k \, \bg_k) > f_{\text{ref}} - c_{\rm{ls}} \, \nu_k \, \|\bg_k\|^2$ \ {\bf do} \ $\nu_k = \sigma_{\rm{ls}} \, \nu_k$ \ {\bf end} \\
{\footnotesize 4}: & \phantom{M} Set \ $\bs_k = -\nu_k \, \bg_k$ \ and update \ $\bx_{k+1} = \bx_k + \bs_k$ \\
{\footnotesize 5}: & \phantom{M} Compute the gradient \ $\bg_{k+1} = \nabla f(\bx_{k+1})$ \\
{\footnotesize 6}: & \phantom{M} {\bf if} \ $\|\bg_{k+1}\| \le {\sf tol} \cdot \|\bg_0\|$, \ {\bf return}, \ {\bf end} \\
{\footnotesize 7}: & \phantom{M} $\by_k = \bg_{k+1} - \bg_k$ \\
{\footnotesize 8}: & \phantom{M} {\bf if} \ $\bs_k^T\by_k < 0$, set $\beta_{k+1} = \max(\min(\|\bg_{k+1}\|_2^{-1}, \, 10^{5}), \ 1)$\\
{\footnotesize 9}: & \phantom{M} {\bf else} \ compute the homogeneous stepsize $\beta_{k+1}$ \eqref{eq:homobb} \ \bf end\\
{\footnotesize 10}: & \phantom{M} Set $\beta_{k+1} = \min(\max(\beta_{k+1}, \, \beta_{\min}), \ \beta_{\max})$\\
\end{tabular}
\end{algorithm}

This algorithm is similar to the one proposed in \cite{daniela2018steplength,raydan1997barzilai}, with the homogeneous stepsize \eqref{eq:homobb} as the key difference in Line~9. In presence of an uphill direction, i.e., when $\bs_k^T\by_k < 0$, the homogeneous stepsize is negative (cf.~Proposition~\ref{prop:prop} (vi)), and therefore must be replaced with a positive quantity. We choose the new stepsize as in \cite{raydan1997barzilai} (cf.~Line~8 of Algorithm~\ref{algo:nonlin}). The convergence of the method is not affected by these choices, since $\beta_k$ stays uniformly bounded, i.e., $\beta_k\in [\beta_{\min}, \beta_{\max}]$ for all $k$. Therefore, the proof of global convergence of Algorithm~\ref{algo:nonlin} can be easily adapted from \cite[Thm.~2.1]{raydan1997barzilai}. While the convergence is not affected, choosing the homogeneous stepsize as starting steplength in the nonmonotone line search might lead to a smaller number of backtracking steps, compared to classical BB stepsizes.

We finally remark that, as for the BB stepsizes, no line search is required for HBB steps when the function $f$ is strictly convex quadratic, i.e., $f(\bx) = \tfrac12\, \bx^TA\bx - \bb^T\bx$, with $A$ SPD (see the results in \cite{dai2002r,raydan1993barzilai}). In fact, since $H_k \equiv A$ is positive definite, from Proposition~\ref{prop:prop} it holds that $\alpha_k \in [\alpha_k^{\rm{BB1}},\,\alpha_k^{\rm{BB2}}]$. This is sufficient to guarantee the R-linear convergence of the corresponding gradient method (cf. \cite[Thm.~4.1]{dai2003alternate} or \cite[Thm.~16]{FHK23}). 

\subsection{Numerical experiments}
\label{sec:exp}
This subsection is devoted to testing the use of the homogeneous stepsize HBB, and an adaptive variant in Algorithm~\ref{algo:nonlin}, on a set of unconstrained optimization problems. We take some general continuously differentiable functions and the suggested starting points therein from the collection in \cite{raydan1997barzilai,andrei2008unconstrained,more1981testing}, as listed in Table~\ref{tab:unconstrained}.
\begin{table}[ht]
\centering
\caption{Unconstrained optimization test problems.}
{\footnotesize
\begin{tabular}{lclc} 
\toprule 
Name & Reference & Name & Reference \\ 
\midrule
{\sf Diagonal 1} & \cite{andrei2008unconstrained} & {\sf Full Hessian FH1} & \cite{andrei2008unconstrained} \\ 
 {\sf Diagonal 2} & \cite{andrei2008unconstrained} & {\sf Full Hessian FH2} & \cite{andrei2008unconstrained} \\ 
 {\sf Diagonal 3} & \cite{andrei2008unconstrained} & {\sf Generalized Rosenbrock} & \cite{andrei2008unconstrained} \\ 
 {\sf Extended Beale} & \cite{andrei2008unconstrained} & {\sf Generalized White and Holst} & \cite{andrei2008unconstrained} \\ 
 {\sf Extended Powell} & \cite{more1981testing} & {\sf Hager} & \cite{andrei2008unconstrained} \\ 
 {\sf Extended Rosenbrock} & \cite{more1981testing} & {\sf Perturbed quadratic} & \cite{andrei2008unconstrained} \\ 
 {\sf Extended White and Holst} & \cite{andrei2008unconstrained} & {\sf Strictly Convex 2} & \cite{raydan1997barzilai} \\ 
\bottomrule
\end{tabular}}
\label{tab:unconstrained}
\end{table}
For all the test functions, the size can be $n\in \{10^2, 10^3, 10^4\}$. The {\sf generalized Rosenbrock}, {\sf generalized White and Holst} and {\sf extended Powell} objective functions have been scaled by the Euclidean norm of the first gradient. 

As for the parameters of Algorithm~\ref{algo:nonlin}, we set the choices $\beta_{\min} = 10^{-30}$, $\beta_{\max} = 10^{30}$, $c_{\rm{ls}} = 10^{-4}$, $\sigma_{\rm{ls}} = \frac12$, $M = 10$, and $\beta_0 = 1$.
The algorithm stops when $\|\bg_k\| \le {\sf tol} \cdot \|\bg_0\|$, with ${\sf tol} = 10^{-6}$, or when $5\cdot 10^4$ iterations are reached. All different steps in Table~\ref{tab:exp} are tested, meaning that we change the stepsize choice in Line~9 of Algorithm~\ref{algo:nonlin} and obtain several gradient methods. Along with the homogeneous stepsize (HBB), BB1 and BB2, we also consider the adaptive ABB method \cite{frassoldati2008new}: 
\begin{equation}
\label{eq:abb}
\beta_k^{\rm ABB} = \left\{ 
\begin{array}{ll}
\min\{\beta_j^{\rm BB2}\mid j = \max\{1,k-m\},\dots,k\}, & \quad {\rm if}\; \beta_k^{\rm BB2} < \eta \, \beta_k^{\rm BB1}, \\[2mm]
\beta_k^{\rm BB1} & \quad \text{otherwise}.
\end{array}
\right.
\end{equation}
Default parameters values are $\eta = 0.8$ and $m = 5$; see, e.g., \cite{daniela2018steplength}. This stepsize is mostly known as $\text{ABB}_{\text{min}}$, while ABB indicates the case for $m = 0$; we choose to indicate it as ABB for any $m$ to ease the notation. Inspired by this, we also propose a straightforward generalization of our HBB step: the adaptive HBB method (AHBB), which takes the stepsize
\begin{equation} \label{ahbb}
\beta_k^{\rm AHBB} = \left\{ 
\begin{array}{ll}
\min\{\beta_j^{\rm HBB}\mid j = \max\{1,k-m\},\dots,k\}, & \ {\rm if}\; \beta_k^{\rm BB2} < \eta \, \beta_k^{\rm BB1}, \\[2mm]
\beta_k^{\rm BB1} & \ \text{otherwise}.
\end{array}
\right.
\end{equation}

\begin{table}[ht]
\centering
\caption{Tested stepsizes.}
\label{tab:exp}
{\footnotesize 
\begin{tabular}{llll} 
\toprule
{\bf Method} & {\bf Description} & {\bf Reference} & \\ 
\midrule
BB1 & Inverse Rayleigh quotient & Cf.~\cite{bb1988} & Eq.~\eqref{bb1} \\
BB2 & Inverse harmonic Rayleigh quotient & Cf.~\cite{bb1988} & Eq.~\eqref{bb1} \\
HBB & Inverse homogeneous Rayleigh quotient && Eq.~\eqref{eq:homobb} \\
ABB & Adaptive: BB1 or (past $m = 5$) BB2 step & Cf.~\cite{frassoldati2008new} & Eq.~\eqref{eq:abb} \\
AHBB & Adaptive: BB1 or (past $m = 5$) HBB step && Eq.~\eqref{ahbb} \\
\bottomrule
\end{tabular}
}
\end{table}
Since the stopping criterion is based on the gradient norm, and $\bs_k = -\nu_k\bg_k$ (cf. Line~4 in Algorithm~\ref{algo:nonlin}), the homogeneous stepsize requires the computation of the two inner products $\by_k^T\by_k$ and $\bg_k^T\by_k$, and therefore it has the same cost as BB2.
In general, all the studied stepsizes in Table~\ref{tab:exp} have very similar cost; BB1 is slightly cheaper, since it does not require the inner product $\by_k^T\by_k$. 

Table~\ref{tab:homounconstr} reports the number of function evaluations and the number of iterations for each problem and stepsize. We remark that, for each problem, all methods converged to the same stationary point. It seems that either HBB or AHBB can be competitive or better than ABB, which is well known for its generally good behavior. For example, HBB behaves nicely in problems {\sf FH1} and {\sf FH2}; AHBB is the best method for the {\sf generalized White and Holst} function and {\sf perturbed quadratic} ($n = 10^4$).

\begin{table}[ht]
\centering
\caption{Number of function evaluations \gf{(NFE)} and iterations for each stepsize and problem. Boldface items are minimal.} 
\label{tab:homounconstr}
\begingroup\footnotesize
\scalebox{0.85}{\begin{tabular}{lr|rrrrr|rrrrr}
 \toprule
 \multicolumn{2}{c}{Problem and size}& \multicolumn{5}{c}{NFE}& \multicolumn{5}{c}{Iterations}\\ \midrule
 && BB1 & BB2 & ABB & HBB & AHBB & BB1 & BB2 & ABB & HBB & AHBB \\ 
 Diagonal 1 &($10^2$) & {\bf 65} & 68 & {\bf 65} & 69 & 67 & {\bf 57} & 63 & 60 & 63 & 62 \\ 
 Diagonal 1 &($10^3$) & 305 & 194 & 160 & {\bf 159} & 165 & 222 & 180 & 149 & {\bf 145} & 156 \\ 
 Diagonal 1 &($10^4$) & 761 & 433 & {\bf 302} & 346 & 310 & 491 & 409 & {\bf 289} & 323 & 297 \\ 
 Diagonal 2 &($10^2$) & 75 & 68 & {\bf 59} & 73 & 66 & 68 & 67 & {\bf 58} & 68 & 63 \\ 
 Diagonal 2 &($10^3$) & 447 & 234 & {\bf 165} & 281 & 212 & 286 & 219 & {\bf 157} & 185 & 175 \\ 
 Diagonal 2 &($10^4$) & 689 & 702 & {\bf 322} & 1221 & 542 & 405 & 663 & \bf{307} & 729 & 407 \\ 
 Diagonal 3 &($10^2$) & 76 & 83 & 73 & {\bf 67} & {\bf 67} & 62 & 73 & 65 & 60 & {\bf 59} \\ 
 Diagonal 3 &($10^3$) & 311 & 200 & 152 & 151 & {\bf 148} & 206 & 184 & 138 & {\bf 137} & {\bf 137} \\ 
 Diagonal 3 &($10^4$) & 349 & 375 & {\bf 238} & 283 & 278 & {\bf 224} & 356 & {\bf 224} & 266 & 263 \\ 
 Ext~Beale &($10^2$) & 50 & 34 & {\bf 33} & {\bf 33} & {\bf 33} & 45 & 30 & 29 & {\bf 27} & 29 \\ 
 Ext~Powell &($10^2$) & 138 & {\bf 80} & 117 & 136 & 156 & 103 & {\bf 74} & 116 & 100 & 131 \\ 
 Ext~Powell &($10^3$) & 117 & {\bf 97} & 114 & 132 & 122 & 97 & {\bf 84} & 113 & 103 & 104 \\ 
 Ext~Powell &($10^4$) & 164 & {\bf 125} & 144 & 161 & 155 & 127 & {\bf 112} & 143 & 114 & 124 \\ 
 Ext~Rosen &($10^2$) & 104 & {\bf 65} & 88 & 69 & 87 & 54 & {\bf 53} & 72 & 55 & 74 \\ 
 Ext~WH &($10^2$) & 99 & 48 & 39 & {\bf 37} & 39 & 62 & 37 & 29 & {\bf 27} & 29 \\ 
 FH1 &($10^2$) & 623 & 304 & 337 & {\bf 256} & 376 & 408 & 279 & 310 & {\bf 221} & 347 \\ 
 FH2 &($10^2$) & 987 & 574 & 563 & {\bf 447} & 499 & 635 & 530 & 537 & {\bf 407} & 470 \\ 
 Gen~Rosen &($10^2$) & 4121 & 3162 & {\bf 2996} & 4808 & 3567 & {\bf 2623} & 2932 & 2736 & 3085 & 2840 \\ 
 Gen~Rosen &($10^3$) & 37769 & {\bf 27324} & 28624 & 37829 & 31516 & {\bf 24020} & 24972 & 26078 & 24082 & 24700 \\ 
 Gen~WH &($10^2$) & 11040 & 8403 & 8866 & 11227 & {\bf 7216} & 6949 & 8107 & 8660 & 7123 & {\bf 5602} \\ 
 Hager &($10^2$) & {\bf 24} & 27 & 27 & 25 & 27 & {\bf 21} & 24 & 24 & 22 & 24 \\ 
 Hager &($10^3$) & {\bf 41} & 46 & 44 & 44 & 45 & {\bf 37} & 42 & 40 & 40 & 41 \\ 
 Hager &($10^4$) & 84 & {\bf 62} & 64 & 66 & 64 & 74 & {\bf 55} & 57 & 59 & 57 \\ 
 Pert quad &($10^2$) & 76 & 98 & {\bf 72} & 76 & 73 & {\bf 63} & 90 & 64 & 66 & 65 \\ 
 Pert quad &($10^3$) & 289 & 220 & 169 & 263 & {\bf 156} & 194 & 198 & 158 & 243 & {\bf 145} \\ 
 Pert quad &($10^4$) & 618 & 350 & 316 & 270 & {\bf 260} & 401 & 330 & 296 & 252 & {\bf 242} \\ 
 S~Conv 2 &($10^2$) & 82 & 62 & {\bf 61} & 66 & 71 & 72 & 58 & {\bf 56} & 57 & 65 \\ 
 S~Conv 2 &($10^3$) & 282 & 156 & {\bf 120} & 156 & 141 & 197 & 145 & {\bf 111} & 137 & 132 \\ 
 S~Conv 2 &($10^4$) & 467 & 234 & {\bf 191} & 213 & 199 & 289 & 221 & {\bf 179} & 195 & 185 \\ 
\bottomrule
\end{tabular}}
\endgroup
\end{table}

\section{Conclusions} \label{sec:con}
We have introduced a new Rayleigh quotient, the homogeneous Rayleigh quotient $\alpha$, which minimizes the homogeneous residual quantity \eqref{homo-sec}. We have shown that the value of $\alpha$ lies between the standard and the harmonic Rayleigh quotient, and so does its sensitivity when approximating extreme eigenvalues of an SPD matrix $A$. Regarding the other eigenvalues and the indefinite case, we have provided examples where the homogeneous Rayleigh quotient may be more accurate than the others.
In addition, in all our experiments, we have noticed that the homogeneous Rayleigh quotient is less sensitive than one of the other two Rayleigh quotients. To some extent, the homogeneous Rayleigh quotient seems to leverage the two Rayleigh quotients: it is less sensitive than the harmonic Rayleigh quotient when estimating the smallest eigenvalues, and less sensitive than the standard Rayleigh quotient when estimating the largest eigenvalues (in magnitude). 

We have derived a {\em nonlinear} Galerkin condition for the homogeneous Rayleigh quotient, in contrast to the linear Galerkin conditions for the standard and the harmonic Rayleigh quotient.
All the results extend to the homogeneous Rayleigh quotient for the generalized eigenvalue problem $(A,B)$, when $B$ is SPD and $A$ is symmetric.

Finally, we have considered the homogeneous Rayleigh quotient as inverse stepsize (HBB) in gradient methods for unconstrained optimization problems.
\gf{This stepsize has first and independently been obtained recently from a different angle by Li, Zhang, and Xia \cite{LZX22, LX21}.}
We have also proposed the AHBB steplength as an alternative to the ABB stepsize, based on the homogeneous stepsize. Experiments show that this variant sometimes performs better than the classical BB steplengths and ABB.

\begin{acknowledgement} \rm
We are grateful to the referees for their very useful suggestions which considerably improved the quality of the paper.
This work has received funding from the European Union's Horizon 2020 research and innovation programme under the Marie Sklodowska-Curie grant agreement No 812912.
\end{acknowledgement}


\end{document}